\documentclass[12pt,reqno]{amsart}
\usepackage{amsmath, amsfonts, amssymb, amsthm, hyperref}
\usepackage{bm}
\usepackage{rotating}
\allowdisplaybreaks[4]
\def\l{\left}
\def\r{\right}
\def\bg{\bigg}
\def\({\bg(}
\def\){\bg)}
\def\t{\text}
\def\f{\frac}

\def\per{{\rm per}}

\def\ls{\le}

\def\ve{\varepsilon}

\def\eq{\equiv}

\def\da{\delta}

\def\Proof{\noindent{\it Proof}}

\def\Z{\mathbb Z}

\def\Q{\mathbb Q}

\def\1{{\bf 1}}

\def\jacob #1#2{\genfrac{(}{)}{}{}{#1}{#2}}
\def\pmod #1{\ ({\rm{mod}}\ #1)}

\def\<{\langle}
\def\>{\rangle}
\theoremstyle{plain}
\newtheorem{theorem}{Theorem}[section]
\newtheorem{lemma}{Lemma}

\theoremstyle{definition}

\theoremstyle{remark}
\newtheorem{remark}{Remark}

\makeatletter
\@namedef{subjclassname@2020}{%
  \textup{2020} Mathematics Subject Classification}
\makeatother
 \vspace{4mm}

\begin{document}
\hbox{Preprint}
\medskip

\title[Some new results on determinants and permanents]
{Some new results on determinants and permanents}

\author{Bo Jiang}
\address {(Bo Jiang) Department of Mathematics, Nanjing
University, Nanjing 210093, People's Republic of China}
\email{bjiang@smail.nju.edu.cn}

\author{Zhi-Wei Sun}
\address {(Zhi-Wei Sun, corresponding author) School of Mathematics, Nanjing
	University, Nanjing 210093, People's Republic of China}
\email{zwsun@nju.edu.cn}

\keywords{Determinant, permanent, Legendre symbol, congruence.
\newline \indent Supported by the National Natural Science Foundation of China (grant no. 11971222).}
\begin{abstract}
In this paper we confirm several conjectures on determinants and permanents.
 For example, we prove that for any prime $p\eq 3\pmod 4$ the number $2\det[a_{jk}]_{0\le j,k\le (p-1)/2}$ is congruent to a square modulo $p$,
where $a_{jk}=\jacob{j+k}{p}+\jacob{j^2+k^2}{p}$ with $\jacob{\cdot}{p}$ the Legendre symbol.
We also prove that $\per[j^{k-1}]_{1\leq j,k\leq n-1}\eq0\pmod n$ for any integer $n>1$ with $n\not\eq2\pmod 4$.
\end{abstract}
\maketitle

\section{Introduction}
\setcounter{lemma}{0}
\setcounter{theorem}{0}
\setcounter{equation}{0}
\setcounter{conjecture}{0}
\setcounter{remark}{0}
\setcounter{corollary}{0}

For a matrix $A=[a_{jk}]_{1\ls j,k\ls n}$ over a field, let $\det(A)$ or $\det[a_{jk}]_{1\ls j,k\ls n}$
denote its determinant.

Let $p$ be an odd prime, and let $(\f{\cdot}p)$ be the Legendre symbol.
If $p\eq1\pmod4$, let $\ve_p$ and $h(p)$ be the fundamental unit and the class number of the real
quadratic field $\Q(\sqrt p)$ respectively, and write
$$\ve_p^{h(p)}=a_p+b_p\sqrt p\ \ \t{with}\ 2a_p,2b_p\in\Z.$$
In 2004 R. Chapman \cite{Ch} proved that
$$\det\l[\l(\f{j+k}p\r)\r]_{0\le j,k\le(p-1)/2}
=\begin{cases}(-1)^{(p+3)/4}2^{(p-1)/2}a_p&\t{if}\ p\eq1\pmod4,
\\2^{(p-1)/2}&\t{if}\ p>3\ \t{and}\ p\eq3\pmod4,
\end{cases}$$
In 2019, Z.-W. Sun \cite{S19} introduced the determinant
$$S_p=\det\l[\l(\f{j^2+k^2}p\r)\r]_{1\le j,k\le (p-1)/2}
\ \ \t{and}\ \ T_p=\det\l[\l(\f{j^2+k^2}p\r)\r]_{0\le j,k\le (p-1)/2},$$
and proved that
$$S_p=\f2{p-1}T_p\ \ \t{and}\ \ \l(\f{T_p}p\r)=\l(\f 2p\r).$$

 The second author \cite[Remark 1.1]{Sun-R} conjectured that $$\det\l[x+\l(\f{j^2+k^2}p\r)+\l(\f{j^2-k^2}p\r)\r]_{1\ls j,k\ls (p-1)/2}
=\l(\f{p-1}2x-1\r)p^{(p-3)/4}
$$ for any prime $p\eq3\pmod4$. This was confirmed by J. Li and H.-L. Wu \cite{LW}.

Our first theorem essentially confirms Conjecture 4.10(ii) of Z.-W. Sun \cite{Sun-R}.

\begin{theorem}\label{Thm1} Let $p$ be a prime with $p\eq\eq 3\pmod 4$. For the matrix
$$A=\l[\l(\f{j+k}p\r)+\l(\f{j^2+k^2}p\r)\r]_{0\le j,k\le(p-1)/2},$$
we have
$$\l(\f{2\det(A)}p\r)\not=-1.$$
\end{theorem}
\begin{remark} Actually, Sun \cite{Sun-R} conjectured $2\det(A)$ is a quadratic residue modulo $p$
but we are unable to prove that $p\nmid \det(A)$.
\end{remark}

For a matrix $A=[a_{j,k}]_{1\le j,k\le n}$ over a field, its {\it permanent} is given by
$$\per(A)=\per[a_{jk}]_{1\le j,k\le n}=\sum_{\sigma\in S_n}\prod_{j=1}^n a_{j,\sigma(j)}.$$

Our second theorem confirms Conjecture 4.2(iii) of Z.-W. Sun \cite{S21} on permanents.

\begin{theorem}\label{Thm2} For any prime $p>3$ and integer $a\not\eq0\pmod p$, we have
$$\sum_{\tau\in S_p\atop p\mid(aj+\tau(j))\ \t{for no}\ j}\prod_{j=1}^p\f1{aj+\tau(j)}\eq0\pmod{p^2}.$$
\end{theorem}

In 2021 the second author \cite{JAC} proved that
$$\per[j^{k-1}]_{1\le j,k\le n}\eq0\pmod n$$
for every integer $n>2$.

Our third theorem confirms Conjecture 4.2(i) of Sun \cite{JAC}.

\begin{theorem}\label{Thm3} Let $n>1$ be an integer with $n\not\equiv 2\pmod4$.
Then we have
	\begin{equation}
		\per[j^{k-1}]_{1\leq j,k\leq n-1}\eq0\pmod n.
	\end{equation}
\end{theorem}

\section{Proof of Theorem $\ref{Thm1}$}
\setcounter{lemma}{0}
\setcounter{theorem}{0}
\setcounter{equation}{0}
\setcounter{conjecture}{0}
\setcounter{remark}{0}
\setcounter{corollary}{0}

\medskip
\noindent {\it Proof of Theorem \ref{Thm1}}. For convenience, we set $n=(p-1)/2$. 
Define
	\begin{equation*} f(x,y)=(x+y)^n+\sum_{l=0}^{\frac{n-1}{2}}\binom{n}{l}x^{2l}
r_{2n-2l}(y)+\sum_{l=\frac{n+1}{2}}^{n}\binom{n}{l}r_{2l}(x)y^{2n-2l},
\end{equation*}
where $r_m(x)$ is the residue of $x^m$ modulo $\prod_{k=0}^n(x-k)$.
For any $0\leq j,k \leq n$, we have
\[\begin{aligned}
	f(j,k)=\ &(j+k)^n+\sum_{l=0}^{(n-1)/{2}}\binom{n}{l}j^{2l}k^{2n-2l}+\sum_{l=(n+1)/{2}}^{n}\binom{n}{l}j^{2l}k^{2n-2l}
	\\=\ &(j+k)^n+(j^2+k^2)^n
\\\eq\ & \jacob{j+k}{p}+\jacob{j^2+k^2}{p} \pmod p .
\end{aligned}\]
It follows that $\det(A)\eq \det(F)\pmod p$, where $F=[f(j,k)]_{0\ls j,k\ls n}$.

Observe that
\[
f(x,y)=\sum_{u,v=0}^{n}b(u,v)x^uy^v
\]
with
\[
b(u,v)=\begin{cases}
0&\t{if}\ 2\nmid u\ \t{and}\ 2\nmid v, \\
\binom{n}{\frac{u}{2}}[x^v]r_{2n-u}(x)+\binom{n}{\frac{v}{2}}[x^u]r_{2n-v}(x)&\t{if}\ 2\mid u\  \t{and}\ 2\mid v,
\\
\da_{u+v,n}\binom{n}{v}+\binom{n}{\frac{u}{2}}[x^v]r_{2n-u}(x)&\t{if}\ 2\mid u\ \t{and}\  2\nmid v,
\\
\da_{u+v,n}\binom{n}{u}+\binom{n}{\frac{v}{2}}[x^u]r_{2n-v}(x)&\t{if}\ 2\nmid u\ \t{and}\ 2\mid v,
\end{cases}
\]
where $\da_{m,n}$ is the Kronecker delta which takes 1 or 0 according as $m=n$ or not,
and $[x^m]P(x)$ is the coefficient of $x^m$ in the polynomial $P(x)$.

Let $B=[b(u,v)]_{0\ls u,v\ls n}$ and $H=[j^u]_{0\ls j,u\ls n}$. Then $F=HBH^T$
and hence $\det(F)=\det(H)^2\det(B)$.

Let $C_0$ denote the matrix $[c_{uv}]_{0\leq u,v\leq n}$, where
\[
c_{uv}=\begin{cases}
	b(2u,2v)&\t{if}\ 0\leq u,v\leq \frac{n-1}{2},
	\\
	b(2u,2v-n)&\t{if}\ 0\leq u\leq \frac{n-1}{2}\ \t{and}\ \frac{n+1}{2}\leq v\leq n,
	\\
	b(2u-n,2v)&\t{if}\ \frac{n+1}{2}\leq u\leq n\ \t{and}\ 0\leq v\leq \frac{n-1}{2},
	\\
	b(2u-n,2v-n)&\t{if}\ \frac{n+1}{2}\leq u,v\leq n.
\end{cases}
\]
 Note that $C_0$ can be obtained via rearranging the rows and columns of $B$, namely, for the $2n\times 2n$ matrix $B$, first move all odd-numbered rows of $B$ sequentially to the first $n$ rows, and then move all odd-numbered columns sequentially to the first $n$ columns. Thus $\det(B)\eq \det(C_0)\pmod p$.

 Observe that
$$C_0=\begin{bmatrix}
	C_1&C_2\\ C_3&C_4
\end{bmatrix}, $$
where
 \begin{gather*}C_1=[c_{uv}]_{0\leq u,v\leq \frac{n-1}{2}}, \  C_2=[c_{uv}]_{0\leq u\leq (n-1)/{2}\atop (n+1)/{2}\leq v\leq n},
 \ C_3=[c_{uv}]_{(n+1)/{2}\leq u\leq n
 \atop 0\leq v\leq (n-1)/2},\  C_4=[c_{uv}]_{\frac{n+1}{2}\leq u,v\leq n}.
  \end{gather*}

Clearly, $C_4=0$ since $2\nmid n$. When $0\leq j\leq \frac{n-1}{2}$ and $\frac{n+1}{2}\leq k\leq n$, we have $$c_{jk}=b(2j,2k-n)=b(2k-n,2j)=c_{kj}.$$ Thus $C_2=C_3^T$ and $$\det(C_0)=(-1)^{(\frac{n+1}{2})^2}\det(C_2)^2=(-1)^{\frac{p+1}{4}}\det(C_2)^2.$$

Since $p\eq 3\pmod 4$, we have $\jacob{-1}{p}=-1$ and $\frac{p^2-1}{8}\eq\frac{p+1}{4} \pmod 2$.
As $\det(A)\eq \det(H)^2(-1)^{\frac{p+1}{4}}\det(C_2)^2 \pmod p$, we obtain
\[\begin{aligned}
	\jacob{2\det(A)}{p}=\ &\jacob{2}{p}\jacob{\det(H)\det(C_2)}{p}^2\jacob{(-1)^{\frac{p+1}{4}}}{p}
	\\=\ &(-1)^{\frac{p^2-1}{8}}(-1)^{\frac{p+1}{4}}\jacob{\det(H)\det(C_2)}{p}^2
	\\=\
	&\jacob{\det(H)\det(C_2)}{p}^2=0\ \t{or}\  1.
\end{aligned}\]
 This concludes our proof of Theorem \ref{Thm1}.
\qed

\section{Proof of Theorem $\ref{Thm2}$}
\setcounter{lemma}{0}
\setcounter{theorem}{0}
\setcounter{equation}{0}
\setcounter{conjecture}{0}
\setcounter{remark}{0}
\setcounter{corollary}{0}
\begin{lemma}\label{Lem3.1}
	Let $p$ be an odd prime, and let $x$ be an integer not divisible by $p$. Then we have
	\begin{equation}
		x^{-1}\eq 2x^{p-2}-x^{2p-3}\pmod{p^2}.
	\end{equation}
\end{lemma}
\Proof. By Fermat's little theorem, $x^{p-1}-1\eq 0\pmod p$. Thus
$$x^{-1}\eq x^{-1}(1-(x^{p-1}-1)^2)\eq 2x^{p-2}-x^{2p-3}\pmod{p^2}. $$
This ends the proof. \qed

\begin{lemma}\label{Lem3.2}
	 Let $p>3$ be a prime. Let $F(x_1,x_2,\dots,x_{p-1})$ be a homogeneous symmetric polynomial 
with integer coefficients. When $\deg(F)\not\eq 0\pmod{p-1}$, we have
	\begin{equation}\label{F1}
		F(1,2,\dots,p-1)\eq 0\pmod{p}.
	\end{equation}
	Moreover, if $\deg(F)\eq -1\pmod{p-1}$ then
	\begin{equation}\label{F2}
		F(1,2,\dots, p-1)\eq 0\pmod{p^2}.
	\end{equation}
\end{lemma}

\Proof. Let $g$ be a primitive root modulo $p$. Since $\gcd(g,p)=1$ and $F(x_1,x_2,\dots,x_{p-1})$ is a homogeneous symmetric polynomial, we have
$$F(1,2,\dots,p-1)\eq F(g,2g,\dots,(p-1)g)\eq g^{\deg(F)}F(1,2,\dots,p-1)\pmod p. $$
If $\deg(F)\not\eq 0\pmod{p-1}$, then $g^{\deg(F)}\not\eq 1\pmod{p}$, and hence \eqref{F1} holds.

Now assume that $\deg(F)\eq -1\pmod{p-1}$. Then $$G(x_1,x_2,\dots,x_p)=F(x_p-x_1,x_p-x_2,\dots,x_p-x_{p-1})$$ is a homogeneous polynomial with $\deg(G)=\deg(F)$, and a symmetric polynomial
in $x_1,x_2,\dots,x_{p-1}$ when $x_p$ is fixed.
Write
$$
	G(x_1,x_2,\dots,x_p)=\sum_{l=0}^{\deg(F)}G_l(x_1,x_2,\dots,x_{p-1})x_p^l,
$$
where $G_l$ is a homogeneous symmetric polynomial in $x_1,x_2,\dots,x_{p-1}$. Since $\deg(F)$ is odd, we obtain $G_0=-F$. As $\deg(G_1)=\deg(F)-1\not\eq 0\pmod{p-1}$ and $G_0=-F$, we have
$$F(p-1,p-2,\dots,1)=G(1,2,\dots,p)\eq G_0(1,2,\dots,p-1)\eq -F(1,2,\dots,p-1)\pmod{p^2}, $$
and hence \eqref{F2} follows.
\qed

\begin{lemma}\label{Lem3.3}
	Let $p>3$ be a prime. Let $F(x_1,x_2,\dots,x_p)$ be a homogeneous symmetric polynomial with integer coefficients. When $\deg(F)\not\eq 0\pmod{p-1}$, we have
	\begin{equation}
		F(1,2,\dots, p)\eq 0\pmod{p}.
	\end{equation}
	Moreover, if $\deg(F)\eq -1\pmod{(p-1)}$ then
	\begin{equation}
		F(1,2,\dots, p)\eq 0\pmod{p^2}.
	\end{equation}
\end{lemma}
\Proof. Write
$$
F(x_1,x_2,\dots,x_p)=\sum_{l=0}^{\deg(F)}F_l(x_1,x_2,\dots,x_{p-1})x_p^l,
$$
where $F_l$ is a homogeneous symmetric polynomial with integer coefficients. It is easy to see that and $\deg(F_0)=\deg(F)$ and $\deg(F_1)=\deg(F)-1$. According to Lemma \ref{Lem3.2}, we have
$$F(1,2,\dots,p)\eq F_0(1,2,\dots,p-1)\eq 0\pmod p$$
when $\deg(F)\not\eq 0\pmod{(p-1)}$,
and
$$F(1,2,\dots,p)\eq F_1(1,2,\dots,p-1)p+ F_0(1,2,\dots,p-1)\eq 0\pmod {p^2}$$
when $\deg(F)\eq -1\pmod{p-1}$. This ends the proof. \qed
\medskip

\noindent {\it Proof of Theorem \ref{Thm2}}.
Let $Q(x)=2x^{p-2}-x^{2p-3}$. Since $p>3$, we have
\begin{equation*}\label{ST} Q(x)\eq\begin{cases}x^{-1}\ \pmod{p^2}&\t{if}\ p\nmid x,\\0\ \ \ \ \pmod{p^2}&\t{if}\ p\mid x.
	\end{cases}.
\end{equation*}
Let $M_{jk}=Q(aj+k)$ and $M=[M_{jk}]_{1\leq j,k\leq p}$. Then
\begin{equation*}
	\sum_{\tau\in S_p\atop p\mid aj+\tau(j)\ \t{for no}\ j}\prod_{j=1}^{p}\frac{1}{aj+\tau(j)}\eq\per(M)\pmod{p^2}.
\end{equation*}

Let $h(x_1,x_2,\dots,x_p)=\per[Q(ax_j+x_k)]_{1\leq j,k\leq p}$. This is a symmetric polynomial in $x_1,x_2,\dots,x_p$. The degree of any term in $Q(ax_j+x_k)$ is either $p-2$ or $2p-3$. Thus, the degree of any term in $h(x_1,x_2,\dots,x_p)$ is congruent to $-1$ modulo $p-1$. Write
$$h(x_1,x_2,\dots,x_p)=\sum_{l=0}^{L}h_l(x_1,x_2,\dots,x_p),$$
where each $h_l$ is a homogeneous symmetric polynomial with $\deg(h_l)\eq -1\pmod{p-1}$.
 By Lemma \ref{Lem3.3}, we have $h_l(1,2,\dots,p)\eq 0\pmod{p^2}$ for any $l=0,\ldots,L$. Thus, $\per(M)\eq h(1,2,\dots,p)\eq 0\pmod{p^2}$
 as desired. \qed

\section{Proof of Theorem $\ref{Thm3}$}
\setcounter{lemma}{0}
\setcounter{theorem}{0}
\setcounter{equation}{0}
\setcounter{conjecture}{0}
\setcounter{remark}{0}
\setcounter{corollary}{0}

\noindent {\it Proof of Theorem} \ref{Thm3}

We distinguish two cases.

{\it Case} 1. $n\eq 0\  \text{or}\ 3\pmod 4$.

Define
$$\tau=\begin{cases}(1,n-1)(2,n-2)\dots(n/2-1,n/2+1)&\t{if}\ 
 n\eq 0\pmod 4,
  \\(1,n-1)(2,n-2)\dots((n-1)/{2},(n+1)/{2})&\t{if}\ n\eq 3\pmod 4.
   \end{cases}$$
   Then $\tau$ is an odd permutation in $S_{n-1}$,  and for $1\leq j\leq n-1$ we have $\tau(j)=n-j$. Let $A_{n-1}$ be the alternating group consisting of all even permutations in $S_{n-1}$. Then $A_{n-1}\tau=\{\sigma\tau:\ \sigma\in A_{n-1}\}$ is the set of all odd permutations in $S_{n-1}$,
   and we have $S_{n-1}=A_{n-1}\bigcup A_{n-1}\tau$.

 Since $n\eq 0\  \text{or}\ 3\pmod 4$, for any $\sigma\in S_{n-1}$, we have
 $$\sum_{j=1}^{n-1}(\sigma(j)-1)=\sum_{j=1}^{n-1}(j-1)=\frac{(n-1)(n-2)}{2}\eq 1\pmod 2.$$

Thus, 
 \[\begin{aligned}
 	\per[j^{k-1}]_{1\leq j,k\leq n-1}=\ &\sum_{\sigma\in S_{n-1}} \prod_{j=1}^{n-1}j^{\sigma(j)-1}
 	\\=\ &\sum_{\sigma\in A_{n-1}} \(\prod_{j=1}^{n-1}j^{\sigma(j)-1}+\prod_{j=1}^{n-1}j^{\sigma\tau(j)-1}\)
 	\\=\
 	&\sum_{\sigma\in A_{n-1}} \(\prod_{j=1}^{n-1}j^{\sigma(j)-1}+\prod_{j=1}^{n-1}(n-j)^{\sigma(j)-1}\)
 	\\ \eq\
 	&\sum_{\sigma\in A_{n-1}} \l(\prod_{j=1}^{n-1}j^{\sigma(j)-1}+(-1)^{\sum_{j=1}^{n-1}(\sigma(j)-1)}\prod_{j=1}^{n-1}j^{\sigma(j)-1}\r)
 	\\\eq \
 	&\sum_{\sigma\in A_{n-1}} \l(\prod_{j=1}^{n-1}j^{\sigma(j)-1}-\prod_{j=1}^{n-1}j^{\sigma(j)-1}\r)\eq 0\pmod{n}.
 \end{aligned}\]

 {\it Case} 2. $n$ is a prime congruent to $1$ modulo $4$, or $n=9$. 
 
 Let $g$ be a primitive root modulo $n$. Since $\gcd(g,n)=1$, we obtain
$$\per[j^{k-1}]_{1\leq j,k\leq n-1}\eq\per[(gj)^{k-1}]_{1\leq j,k\leq n-1}\eq g^{\frac{(n-1)(n-2)}{2}}\per[j^{k-1}]_{1\leq j,k\leq n-1}\pmod n.$$

Let $\varphi$ be Euler's totient function. When $n$ is prime, we have $\varphi(n)=n-1$ and $$\frac{(n-1)(n-2)}{2}\eq\frac{n-1}{2}\not\eq 0\pmod{\varphi(n)}.$$ When $n=9$, we have $\varphi(n)=6$ and $\frac{(n-1)(n-2)}{2}\eq4\not\eq 0\pmod{\varphi(n)}$. Thus, we have $g^{\frac{(n-1)(n-2)}{2}}\not\eq 1\pmod n$, and hence
$\per[j^{k-1}]_{1\leq j,k\leq n-1}\eq 0\pmod n.$
\medskip

{\it Case} 3. $n\eq1\pmod4$, but $n$ is neither a prime nor equal to $9$.

Let $p$ be any prime divisor of $n$, and write $n=p^{v_p(n)}n_p$,
 where $v_p(n)$ is the $p$-adic valuation of $n$, and $n_p$ a positive integer not divisible by $p$.

When $n_p>1$, for any $\sigma\in S_{n-1}$, we have
$$v_p\(\prod_{j=1}^{n-1}j^{\sigma(j)-1}\)\geq \sum_{l=1}^{v_p(n)}v_p(p^l)(\sigma(p^l)-1)+v_p(2p)(\sigma(2p)-1)\geq v_p(n). $$

 When $n_p=1$, if $p>3$ and $v_p(n)\geq 2$, or $p=3$ and $v_p(n)\geq 3$, then for any $\sigma\in S_{n-1}$ we have
$$v_p\(\prod_{j=1}^{n-1}j^{\sigma(j)-1}\)\geq \sum_{l=1}^{v_p(n)-1}v_p(p^l)(\sigma(p^l)-1)+v_p(2p)(\sigma(2p)-1)++v_p(4p)(\sigma(4p)-1)\geq v_p(n). $$

 As $p^{v_p(n)}\mid \per[i^{j-1}]_{1\leq i,j\leq n-1}$ for any prime divisor $p$ of $n$, we have $n\mid \per[i^{j-1}]_{1\leq i,j\leq n-1}$ as desired.
\qed

\end{document}